\documentclass[a4paper,12pt]{amsart}
\usepackage{ifthen}
\usepackage{graphicx}
\usepackage{mathrsfs}
\nonstopmode
\numberwithin{equation}{section}
\newtheorem{thm}{Theorem}[section]
\newtheorem{cor}[thm]{Corollary}
\newtheorem{lem}[thm]{Lemma}

\newtheorem{Cor}{Corollary}

\theoremstyle{definition}

\newenvironment{pf}[1][]{%
 \vskip 3mm
 \noindent
 \ifthenelse{\equal{#1}{}}%
  {{\slshape Proof. }}%
  {{\slshape #1.} }%
 }%
{\qed\bigskip}

\newcounter{alphabet}
\newcounter{tmp}
\newenvironment{Thm}[1][]{\refstepcounter{alphabet}%
\bigskip%
\noindent%
{\bf Theorem \Alph{alphabet}}%
\ifthenelse{\equal{#1}{}}{}{ (#1)}%
{\bf .}
\itshape}{\vskip 8pt}

\newcommand{\A}{{\mathcal A}}

\newcommand{\C}{{\mathbb C}}
\newcommand{\D}{{\mathbb D}}

\newcommand{\F}{{\mathcal F}}

\newcommand{\K}{{\mathcal K}}

\newcommand{\R}{{\mathbb R}}

\newcommand{\es}{{\mathcal S}}

\newcommand{\CC}{{\mathcal C}}

\newcommand{\bD}{{\overline{\mathbb D}}}

\renewcommand{\Im}{{\,\operatorname{Im}\,}}

\renewcommand{\Re}{{\,\operatorname{Re}\,}}

\renewcommand{\arg}{\,{\operatorname{arg}\,}}

\newcommand{\aand}{{\quad\text{and}\quad}}

\newcounter{minutes}
\divide\time by 60
\newcounter{hours}
\multiply\time by 60
\addtocounter{minutes}{-\time}

\begin{document}
\bibliographystyle{amsplain}
\title{
An extremal problem for univalent functions
}


\author[T.~Sugawa]{Toshiyuki Sugawa}
\address{Graduate School of Information Sciences,
Tohoku University, Aoba-ku, Sendai 980-8579, Japan}
\email{sugawa@math.is.tohoku.ac.jp}
\author[L.-M.~Wang]{Li-Mei Wang}
\address{School of Statistics,
University of International Business and Economics, No.~10, Huixin
Dongjie, Chaoyang District, Beijing 100029, China}
\email{wangmabel@163.com} 
\begin{abstract}

For a real constant $b,$ we give sharp estimates of
$\log|f(z)/z|+b\arg[f(z)/z]$ for subclasses of normalized univalent
functions $f$ on the unit disk.
\end{abstract}
\thanks{
The present research was supported in part by JSPS Grant-in-Aid for
Scientific Research (B) 22340025, and  also by National Natural
Science Foundation of China (No.~11326080)} \maketitle

\section{Introduction}

Let $\A$ denote the class of analytic functions $f$ on the unit
disk $\D=\{z\in\C: |z|<1\}$ normalized so that $f(0)=0$ and $f'(0)=1.$
The subclass $\es$ of $\A$ consisting of all univalent functions
has attracted much interest for many years in the univalent function theory.
In the present paper, we are primarily interested in the extremal problem
to find the value of
$$
\Psi_z(t,\F)=\sup_{f\in\F}\Re\left[e^{it}\log\frac{f(z)}{z}\right]
$$
for a subclass $\F$ of $\es$ and $z\in\D.$
Here and hereafter, $g(z)=\log[f(z)/z]=\log|f(z)/z|+i\arg[f(z)/z]$ will be
understood as the holomorphic branch of logarithm determined by $g(0)=0.$
It is often more convenient to consider the quantities
$$
\Phi_z^+(b,\F)=\sup_{f\in\F}\left\{\log\left|\frac{f(z)}{z}\right|
+b\arg\frac{f(z)}{z}\right\}
$$
and
$$
\Phi_z^-(b,\F)=\inf_{f\in\F}\left\{\log\left|\frac{f(z)}{z}\right|
+b\arg\frac{f(z)}{z}\right\}.
$$
Then, we have the obvious relation
\begin{equation}\label{eq:rel}
\Psi_z(t,\F)=\begin{cases}
(\cos t)\Phi_z^+(-\tan t,\F) &\quad\text{if}~\cos t>0, \\
\null&\null \\
(\cos t)\Phi_z^-(-\tan t,\F) &\quad\text{if}~\cos t<0.
\end{cases}
\end{equation}
Therefore, the first problem is essentially equivalent to finding
the values of $\Phi_z^\pm(b,\F)$ (except for the case when $e^{it}=\pm i$).
We also consider the quantities
$$
\Psi(t,\F)=\sup_{z\in\D}\Psi_z(t,\F)
$$
and
$$
\Phi^+(b,\F)=\sup_{z\in\D}\Phi_z^+(b,\F)
\aand
\Phi^-(b,\F)=\inf_{z\in\D}\Phi_z^-(b,\F).
$$

The above extremal problems will reduce to geometric ones, once we know
about the shape of the variability region $W_z(\F)$
of $\log[f(z)/z]$ for a subclass $\F$ of $\es$ and a fixed point $z\in\D$
defined by
$$
W_z(\F)=\left\{\log\frac{f(z)}{z}: f\in\F\right\}.
$$
Indeed, for instance, we have
$$
\Psi_z(t,\F)
=\sup_{w\in W_z(\F)}\Re[e^{it}w]
=\sup_{u+iv\in W_z(\F)}(u\cos t-v\sin t).
$$

We note that $W_z(\F)=W_r(\F)$ for $r=|z|$ if $\F$ is rotationally
invariant; in other words, if the function
$e^{-i\theta}f(e^{i\theta}z)$ belongs to $\F$ whenever $f\in\F$ and
$\theta\in\R.$ The union
$$
W(\F)=\bigcup_{z\in\D}W_z(\F)
$$
is called the full variability region of $\log[f(z)/z]$ for $\F.$

In the present paper, we will discuss those regions of variability
and the corresponding extremal problems for typical subclasses of $\es.$

\section{Univalent functions}

Grunsky \cite{Gru32} gave a description of $W_z(\es)$
(see also \cite[\S 10.9]{Duren:univ}).

\begin{Thm}[Grunsky]
For $z\in\D$ with $r=|z|,$ the variability region $W_z(\es)$ of $\log[f(z)/z]$
for $\es$ is the closed disk
$$
\left|w-\log\frac1{1-r^2}\right|\le \log\frac{1+r}{1-r}.
$$
\end{Thm}

In particular, we see that for a fixed $z\in\D$ with $|z|=r$
and a real constant $t,$
$$
\Psi_z(t,\es)
=(\cos t)\log\frac1{1-r^2}+\log\frac{1+r}{1-r}
=(1-\cos t)\log(1+r)-(1+\cos t)\log(1-r).
$$
For more general extremal problems, the reader may refer to the monograph
\cite{Alek:para} by Alexsandrov.

In particular, letting $r\to1^-$ in the above, we obtain
$\Psi_z(t,\es)\to+\infty$ if $\cos t>-1$
and $\Psi_z(\pi,\es)=2\log(1+r)\to 2\log2=\log4.$
Hence, in view of \eqref{eq:rel}, we obtain $\Phi^+(b,\es)=+\infty$ for $b\in\R$
and
$$
\Phi^-(b,\es)=\begin{cases}
-\infty & \text{if}~ b\ne0, \\
-\log4 & \text{if}~ b=0.
\end{cases}
$$

More precisely, as a corollary of the Grunsky theorem, we have the following.

\begin{Cor}
The full variability region $W(\es)$ of $\log[f(z)/z]$ for $\es$ is
the half-plane $\{w: \Re w>-\log 4\}.$
\end{Cor}

\begin{pf}
Since $\es$ is rotationally invariant, $W_z(\es)=W_r(\es)$ for
$r=|z|.$ For a fixed $\eta\in\R,$ the intersection of the disk
$W_r(\es)$ with the horizontal line $\Im w=\eta$ is the segment with
the endpoints
$$
-\log(1-r^2)\pm\sqrt{\left(\log\frac{1+r}{1-r}\right)^2-\eta^2}+\eta i
$$
for $r$ so close to $1$ that $|\eta|\le\log[(1+r)/(1-r)].$ Since
$$
-\log(1-r^2)-\sqrt{\left(\log\frac{1+r}{1-r}\right)^2-\eta^2}\to-\log4
$$
whereas
$$
-\log(1-r^2)+\sqrt{\left(\log\frac{1+r}{1-r}\right)^2-\eta^2}\to+\infty
$$
as $r\to1^-,$ we see that $\{w\in W(\es): \Im w=\eta\}=\{x+i\eta: x>-\log4\}.$
The proof is now complete.
\end{pf}

\section{Starlike functions  and convex functions}

For the subclasses $\es^{*}$ and $\K$ of $\es$ consisting of
starlike and convex functions respectively, Marx \cite[Satz B, C]{Marx32}
obtained essentially the following result.

\begin{Thm}[Marx]  For $z\in\D$ with $r=|z|,$ the variability region
$W_z(\es^{*})$ of $\log[f(z)/z]$ for $\es^{*}$ is given as
$$
\left\{-2\log(1-\zeta) :\, |\zeta|\leq r\right\}
$$
and $W_z(\K)$ of $\log[f(z)/z]$ for $\K$ is
$$
\left\{-\log(1-\zeta) :\, |\zeta|\leq
r\right\}=\frac{1}{2}W_z(\es^{*}).
$$
\end{Thm}

Note that $W_r(\es^*)$ is nothing but the image of the disk $|\zeta|\le r$
under the mapping $\log[k(z)/z],$ where $k(z)$ is
the Koebe function; namely, $k(z)=\frac{z}{(1-z)^2}.$
Similarly, $W_r(\K)$ is the image of the disk $|\zeta|\le r$ under the mapping
$\log[l(z)/z],$ where $l(z)=z/(1-z)$ is an extremal convex function.
It is easy to see that the regions $W_r(\es^{*})$ and $W_r(\K)$ are
convex and symmetric with respect to the real axis.

As an application of the Marx theorem, we can solve the corresponding extremal
problem.
We present it only for starlike functions since we have only to take the half for 
convex functions.

\begin{thm}
For a fixed
$z\in\D$ with $|z|=r$ and a real number $b,$
$$
\Phi_z^+(b,\es^*)
=\log(1+b^2)-2\log(\sqrt{\triangle}-r)+2b\arctan\frac{br}{\sqrt{\triangle}}
$$
and
$$
\Phi_z^-(b,\es^*)
=\log(1+b^2)-2\log(\sqrt{\triangle}+r)-2b\arctan\frac{br}{\sqrt{\triangle}},
$$
where $\triangle=1+b^2(1-r^2)$.
\end{thm}

\begin{pf} By theorem B, we have for $z\in \D$ with $|z|=r$,
$$
\Phi_z^+(b,\es^*)
=\sup_{f\in\es^*}\left\{\log\left|\frac{f(z)}{z}\right|
+b\arg\frac{f(z)}{z}\right\}=\sup_{|\xi|\leq
r}\left\{-2\log|1-\xi|-2b\arg(1-\xi)\right\}.
$$
By making use of the maximum principle for harmonic functions,
$$
\Phi_z^+(b,\es^*) =\sup_{|\xi|=
r}\left\{-2\log|1-\xi|-2b\arg(1-\xi)\right\}.
$$
The same argument yields
$$
\Phi_z^-(b,\es^*) =\inf_{|\xi|=
r}\left\{-2\log|1-\xi|-2b\arg(1-\xi)\right\}.
$$

 For $\xi=re^{i\theta}$, let
\begin{eqnarray*}
\phi(\theta)&=&-2\log|1-re^{i\theta}|-2b\arg(1-re^{i\theta})\\
&=&-\log(1+r^2-2r\cos\theta)+2b\arctan\frac{r\sin\theta}{1-r\cos\theta}.
\end{eqnarray*}
We need to find the supremum and infimum of $\phi(\theta)$ over $\R.$
Since $\phi$ is periodic with period $2\pi,$ it is enough to find
(local) maxima and minima in the interval $[0,2\pi)$.
A simple calculation yields
$$
\phi'(\theta)=-2r\cdot\frac{\sin\theta-b\cos\theta+br}{1-2r\cos\theta+r^2}.
$$
Thus extremal values of $\phi(\theta)$ are attained at the points
$\theta$ satisfying
$$
\sin\theta-b\cos\theta+br=0.
$$
By solving the above equation, we have $\theta=\theta_1, \theta_2,$
where
$$
\cos\theta_1=\frac{b^2r+\sqrt{\triangle}}{1+b^2},\quad
\sin\theta_1=\frac{b\sqrt{\triangle}-br}{1+b^2}
$$
and
$$
\cos\theta_2=\frac{b^2r-\sqrt{\triangle}}{1+b^2},\quad
\sin\theta_2=\frac{-b\sqrt{\triangle}-br}{1+b^2}.
$$
We note that such $\theta_1, \theta_2$ exist uniquely on $[0,2\pi)$ since
$$
\left(\frac{b^2r\pm\sqrt{\triangle}}{1+b^2}\right)^2
+\left(\frac{\pm b\sqrt{\triangle}-br}{1+b^2}\right)^2=1.
$$
Thus
 $$
\Phi_z^+(b,\es^*)=\max\{\phi(\theta_1),\phi(\theta_2)\}=\log(1+b^2)-2\log(\sqrt{\triangle}-r)+2b\arctan\frac{br}{\sqrt{\triangle}}
 $$
and
 $$
\Phi_z^-(b,\es^*)=\min\{\phi(\theta_1),\phi(\theta_2)\}=\log(1+b^2)-2\log(\sqrt{\triangle}+r)-2b\arctan\frac{br}{\sqrt{\triangle}}.
 $$
 The proof is now completed.
\end{pf}

We observe that 
$$
\Phi_z^+(b,\es^*)+\Phi_z^-(b,\es^*)=-2\log(1-r^2),
$$
which is independent of the parameter $b.$
Since
$$
\frac{\partial}{\partial b}\Phi_z^+(b,\es^*)
=-\frac{\partial}{\partial b}\Phi_z^-(b,\es^*)
=2\arctan\frac{br}{\sqrt{1+b^2-b^2r^2}},
$$
we can see that $\Phi_z^+(b,\es^*)$ is increasing in $b>0$ and decreasing in $b<0.$
In particular,
$$
\Phi_z^+(b,\es^*)\geq \Phi_z^+(0,\es^*)=-2\log(1-r)
$$
and
$$
\Phi_z^-(b,\es^*)\leq \Phi_z^-(0,\es^*)=-2\log(1+r).
$$

Letting $r\to1^-,$ we obtain the following corollary.

\begin{cor}\label{cor:starlike}
For a real number $b,$
$$
\Phi^+(b,\es^*)=+\infty
\aand
\Phi^-(b,\es^*)=\log(1+b^2)-\log4-2b\arctan b.
$$

\end{cor}

Theorem 3.1 and Corollary 3.2 assure the following result.

\begin{cor}
For a fixed $z\in\D$ with $|z|=r$ and a real number $b$,
$$
\Phi_z^+(b,\K) =\frac{1}{2}\Phi_z^+(b,\es^*),\quad
\Phi_z^-(b,\K)=\frac{1}{2}\Phi_z^-(b,\es^*)
$$
and
$$
\Phi^+(b,\K)=+\infty, \quad
\Phi^-(b,\K)=\frac{1}{2}\log(1+b^2)-\log2-b\arctan b.
$$
\end{cor}

\section{Close-to-convex functions}

Biernacki \cite{Bier36} determined the variability region $W_z(\CC)$
of $\log[f(z)/z]$ for the class $\CC$ of linearly accessible
functions (now known as {\it close-to-convex} functions). That is,
$W_z(\CC)=\{-\log[2u^2/(u+v)]: |u-1|\le |z|, |v-1|\le |z|\}.$ He
also showed that $W(\CC)\subset\{w: |\Im w|<3\pi/2\}.$
Since that is somewhat implicit, Kato and the authors \cite{KSW14} offered
another expression for it; that is, $W_z(\CC)=h(\bD_r)$
for $r=|z|,$ where $\bD_r=\{z: |z|\le r\}$
and $h(z)=\log(1+ze^{2i\phi})-3\log(1+z),~ \phi=\arg(1+z/3).$
It is, however, still difficult to compute $\Phi_z^\pm(b,\CC)$ for $z\in\D.$
Thus, our main concern in the
present paper will be determination of the quantity
$\Phi^\pm(b,\CC)$ because we have a relatively simple expression of
$W(\CC).$

\begin{lem}[Theorem 1.4 in \cite{KSW14}]\label{lem:LU1}
The full variability region $W(\CC)$ for close-to-convex functions
is the unbounded Jordan domain whose boundary is the
Jordan arc $-\gamma((-2\pi,2\pi)).$
Here,
$$
\gamma(t)=
\begin{cases}
\log(1+3e^{it}) & \quad\text{if}~ |t|< \pi \\
\log(1-e^{it})+\dfrac{t}{|t|}\pi i & \quad\text{if}~ \pi\le|t|<2\pi.
\end{cases}
$$
\end{lem}

Note that the region $W(\CC)$ is contained in the parallel strip
$\{w: |\Im w|<3\pi/2\}$ as was already shown by Biernacki \cite{Bier36}.
By making use of the above lemma, we now describe $\Phi^\pm(b,\CC).$

\begin{thm}\label{thm:main}
Let $b$ be a real number.
Then, $\Phi^+(b,\CC)=+\infty$ and
$$
\Phi^-(b,\CC)=\begin{cases}
-\dfrac12\log\dfrac{2(5-4b^2+3\sqrt{1-8b^2})}{1+b^2}
-b\arctan\dfrac{3b}{\sqrt{1-8b^2}}
& \text{if}~ |b|\le b_0, \\
\null & \null \\
\dfrac12\log(1+b^2)-\log2-|b|(\arctan |b|+\pi) & \text{if}~|b|\ge b_0.
\end{cases}
$$
Here, $b_0=0.24001\dots$ is the unique solution to the equation
$$
2b\left(\arctan b-\arctan\dfrac{3b}{\sqrt{1-8b^2}}+\pi\right)
=\log(5-4b^2+3\sqrt{1-8b^2})-\log 2
$$
in $0<b<1/2\sqrt2.$
\end{thm}

\begin{pf}
Since $\es^*\subset\CC,$ Corollary \ref{cor:starlike} yields
$\Phi^+(b,\CC)\ge\Phi^+(b,\es^*)=+\infty$ for every $b.$

We next consider $\Phi^-(b,\CC).$
For brevity, we put $\Phi(b)=\Phi^-(b,\CC)$ throughout the proof.

In order to prove the theorem, we translate our problem
into a geometric one concerning the curve $\gamma$ given in Lemma \ref{lem:LU1}.
First of all, we note that $\gamma$ is symmetric in the sense that
$\gamma(-t)=\overline{\gamma(t)}.$ Therefore, it is enough to
consider the case $0\le t\le2\pi$ unless otherwise stated. We now
study the regularity of the curve $\gamma(t)$ at $t=\pi.$ A direct
computation shows that the left and right tangent vectors
$\gamma'(\pi^-)=3i/2$ and $\gamma'(\pi^+)=i/2$ have the same
direction. Therefore, by a re-parametrization, we see that the
boundary of $W(\CC)$ is of class $C^1.$ We remark, however, that it
is not of class $C^2.$ Indeed, this can be confirmed by observing
that $\exp(\gamma([0,\pi]))$ and $\exp(\gamma([\pi,2\pi]))$ are
(half-)circles with different radii.

We next study convexity of the curve $\gamma.$
In \cite{KSW14}, we already saw that the curve $\gamma$ is not convex.
More precisely, we compute
$$
\frac{d}{dt}\arg \gamma'(t)=
\begin{cases}
\dfrac{1+3\cos t}{|1+3e^{it}|^2} & \text{if}~ 0\le t<\pi, \\
\null & \null \\
\dfrac{1-\cos t}{|1-e^{it}|^2} & \text{if}~\pi<t<2\pi.
\end{cases}
$$
Therefore, the curve $\gamma$ is convex in $0<t<\arccos(-1/3)$ and
$\pi<t<2\pi$ and concave in $\arccos(-1/3)<t<\pi.$ It is important
in the sequel to find the exact form of the convex hull
$\widehat\Omega$ of $\Omega.$ The newly added boundary
$\partial\widehat\Omega-\partial\Omega$ consists of the line segment
joining the two points of tangency of a common tangent line to
$\gamma$ on two parts $0<t<\arccos(-1/3)$ and $\pi<t<2\pi,$ and its
reflection in the real axis.

We should thus find the common tangent line. Let $\gamma(u)$ and
$\gamma(v)$ be the points of tangency of the common tangent line,
where $0<u<\arccos(-1/3)$ and $\pi<v<2\pi.$ Necessary conditions are
described by
\begin{equation}\label{eq:1}
\arg\gamma'(u)=\arg\gamma'(v)=\arg[\gamma(v)-\gamma(u)].
\end{equation}
Since $\gamma'(v)=e^{iv/2}/(2\sin(v/2)),$ we have
\begin{equation}\label{eq:2}
\frac v2=\arg\gamma'(v)=\arg\gamma'(u)=u+\frac\pi2-\alpha,
\end{equation}
where
$$
\alpha=\arg(1+3e^{iu})=\arctan\frac{3\sin u}{1+3\cos u}.
$$
Simple computations give us
$$
\Re[\gamma(v)-\gamma(u)]=\log\left|\frac{1-e^{iv}}{1+3e^{iu}}\right|
=\frac12\log\frac{1-\cos v}{5+3\cos u}
$$
and
$$
\Im[\gamma(v)-\gamma(u)]=\arg(1-e^{iv})+\pi-\arg(1+3e^{iu})
=\frac{v+\pi}2-\alpha.
$$
Hence, the second equation in \eqref{eq:1} yields the relation
\begin{equation}\label{eq:3}
\tan\frac v2=
\frac{v+\pi-2\alpha}{\log(1-\cos v)-\log(5+3\cos u)}.
\end{equation}
In view of \eqref{eq:2}, we have
\begin{equation}\label{eq:4}
\tan\frac v2=\cot(\alpha-u)=\frac{1+\tan\alpha\tan u}{\tan\alpha-\tan u}
=-\frac{3+\cos u}{\sin u}
\end{equation}
and
$$
1-\cos v=2\sin^2\frac v2=2\cos^2(\alpha-u)=\frac{(3+\cos u)^2}{5+3\cos u}.
$$
Substituting these and \eqref{eq:2} into \eqref{eq:3}, we obtain
$$
-\frac{3+\cos u}{\sin u}
=\frac{u+\pi-2\alpha}{\log(3+\cos u)-\log(5+3\cos u)}.
$$
We summarize the above observations as follows.
The slope of the tangent line to $\gamma$ increases from
$-\infty$ to $-(3+\cos u)/\sin u=\tan(v/2)$ as $t$ moves from 0 to $u.$
The tangent line to $\gamma$ at $t=u$ is tangent, at the same time, to
$\gamma$ at $t=v.$
The part $\gamma((u,v))$ is thus contained in the interior of
$\widehat\Omega.$
The slope of the tangent line to $\gamma$ increases from $\tan(v/2)$ to 0
as $t$ moves from $v$ to $2\pi.$

Let $\Omega$ be the domain $\{-w: w\in W(\CC)\}.$
By Lemma \ref{lem:LU1}, $\Omega$ is an unbounded Jordan domain
bounded by the curve $\gamma$ with $0\in\Omega.$
Then,
$$
\Phi(b)=\inf_{w\in W(\CC)}(\Re w+b \Im w)
=-\sup_{X+iY\in\Omega}(X+bY)
=-\max_{X+iY\in\overline\Omega}(X+bY)
$$
Here, we recall that $\Omega$ is contained in the region $X<\log4, |Y|<3\pi/2.$
Hence, the supremum was able to be replaced by the maximum above by taking the
points over the closure of $\Omega.$

Since $\Omega$ is symmetric in the real axis, we have $\Phi(-b)=\Phi(b).$
Hence, we may assume that $b\ge0$ in the proof of Theorem \ref{thm:main}.
When $b=0,$ obviously $\Phi(0)=-\gamma(0)=-\log 4,$ which agrees with
the assertion of the theorem.
In the sequel, we thus assume that $b>0.$
For a given $b>0,$ let $Z_0=X_0+iY_0$ be a point in
$\overline\Omega$ at which $X+bY$ takes its maximum over all
$X+iY\in\overline\Omega.$ It is obvious that $Z_0\in\partial\Omega$
with $Y_0>0$ and that the line $X+bY=X_0+bY_0(=-\Phi(b))$ is tangent
to the curve $\gamma$ at $Z_0.$ Since $Z_0$ is a support point for
the functional $X+bY$ over $\overline\Omega,$ $Z_0=\gamma(t)$ for
some $t$ with $0<t\le u$ or $v\le t<2\pi,$ where $u$ and $v$ are as
above.

When $t\le u(<\pi),$ we have
$X_0=\log|1+3e^{it}|=\frac12\log2(5+3\cos t), Y_0=\arg(1+3e^{it}).$
Also, since the slope of the line $X+bY=-\Phi(b)$ is $-1/b,$ we have
the relation
$$
-\frac1b=\frac{\Im\gamma'(t)}{\Re\gamma'(t)} =\frac{3+\cos t}{-\sin t},
$$
which is
$$
b=\frac{\sin t}{3+\cos t}.
$$
Thus we can get the following equation of $\cos t$
$$
\cos^2 t+ b^2(3+\cos t)^2=1.
$$
Therefore
$$
\cos t=\frac{-3b^2\pm \sqrt{1-8b^2}}{1+b^2}.
$$
Since in this case $\cos t>-1/3$, we have
$$
\cos t=\frac{-3b^2+\sqrt{1-8b^2}}{1+b^2}
$$
and
$$
\tan Y_0=\frac{3\sin t}{1+3\cos t}=\frac{3b(3+\cos t)}{1+3\cos
t}=\frac{3b}{\sqrt{1-8b^2}}.
$$
Hence,
\begin{align*}
\Phi(b)&=-(X_0+bY_0) \\
&=-\dfrac12\log\dfrac{2(5-4b^2+3\sqrt{1-8b^2})}{1+b^2}
-b\arctan\dfrac{3b}{\sqrt{1-8b^2}}=:p(b)
\end{align*}
as is stated in the theorem.

When $t\ge v(>\pi),$ we have $X_0=\log|1-e^{it}|=\log[2\sin(t/2)]$
and $Y_0=\arg(e^{it}-1)=(t+\pi)/2.$ Similarly, we have the relation
$-\frac1b=(1-\cos t)/\sin t=\tan(t/2),$ which is equivalent to
$b=-\cot\frac t2=\tan\frac{t-\pi}2.$ Therefore,
$$
\Phi(b)=-(X_0+bY_0)=-\log\frac2{\sqrt{1+b^2}}-b(\arctan
b+\pi)=:q(b).
$$
Let $b_0=\sin u/(3+\cos u)=-\cot(v/2).$ Then $b_0$ must satisfy the
relation $p(b_0)=q(b_0).$ Indeed, $b_0$ is a unique solution to the
equation $p(b)=q(b)$ in $0<b<1/\sqrt{8}$, since
$$
q'(b)-p'(b)=\dfrac{\left(24b^2-3+(4b^2-5)\sqrt{1-8b^2}\right)(\pi+\arctan
b-\arctan
\frac{3b}{\sqrt{1-8b^2}})}{\sqrt{1-8b^2}(5-4b^2+3\sqrt{1-8b^2})}<0,
$$
$q(0)-p(0)=\log 2>0$ and 
$$
\lim_{b\to(1/\sqrt{8})^{-}}(q(b)-p(b))=\log\frac32-
\frac1{2\sqrt2}\left(\arctan\frac1{2\sqrt2}+\frac\pi2\right)<0.
$$
The proof of Theorem \ref{thm:main} has been completed.
\end{pf}

In view of the relation \eqref{eq:rel}, we obtain the following.

\begin{cor}
For a real constant $t$ with $|t|<\pi/2$,
\begin{eqnarray*}
&&
\Psi(t,\CC) \\
&=&\begin{cases}
-\dfrac{1}{2}(\cos t)\log\big[2(9\cos^2t-4+3\cos^2t\sqrt{1-8\tan^2t})\big]
-(\sin t)\arctan\dfrac{3\tan t}{\sqrt{1-8\tan^2t}} & \null \\
\qquad\qquad\qquad\qquad\qquad\qquad\qquad\qquad\qquad\qquad
\text{if}~ |\tan t|\le b_0, &\null\\
\null & \null \\
-(\cos t)\log[2\cos t]-|\sin t|(|t|+\pi) \qquad\qquad\quad
 \text{if}~|\tan t|\ge b_0, &\null
\end{cases}
\end{eqnarray*}
where $b_0$ is given in Theorem \ref{thm:main}.
\end{cor}

\section{Application to power deformations}

As an application of the main theorem, we consider power deformations
of a univalent function.
Let $c=a+bi$ be a complex number.
The power deformation of a function $f\in\es$ with exponent $c$
is defined by
$$
f_c(z)=z\left(\frac{f(z)}{z}\right)^c
=z\exp(c \log[f(z)/z]).
$$
See \cite{KS11PT} and \cite{KS12PT} for details about the power deformation.
We now have
$$
|f_c(z)|=|z|\exp(a\log|f(z)/z|-b\arg[f(z)/z]).
$$
Therefore, as a corollary of Theorem \ref{thm:main}, we obtain the following.

\begin{thm}
Let $c=a+bi$ be a complex number.
If $a>0$,
$$
\inf_{f\in\CC\atop z\in\D}\log\left|\frac{f_c(z)}{z}\right|=
\begin{cases}
-\dfrac{a}{2}\log\dfrac{2(5a^2-4b^2+3a\sqrt{a^2-8b^2})}{a^2+b^2}-
b\arctan\dfrac{3b}{\sqrt{a^2-8b^2}}
& \text{if}~ \left|\frac{b}{a}\right|\le b_0, \\
\null & \null \\
\dfrac{a}{2}\log(a^2+b^2)-a\log2a^2-|b|(\arctan
\left|\frac{b}{a}\right|+\pi) &
\text{if}~\left|\frac{b}{a}\right|\ge b_0,
\end{cases}
$$
and, if $a<0$,
$$
\sup_{f\in\CC \atop z\in\D}\log\left|\frac{f_c(z)}{z}\right|=
\begin{cases}
-\dfrac{a}{2}\log\dfrac{2(5a^2-4b^2+3a\sqrt{a^2-8b^2})}{a^2+b^2}-
b\arctan\dfrac{3b}{\sqrt{a^2-8b^2}}
& \text{if}~ \left|\frac{b}{a}\right|\le b_0, \\
\null & \null \\
\dfrac{a}{2}\log(a^2+b^2)-a\log2a^2-|b|(\arctan
\left|\frac{b}{a}\right|+\pi) &
\text{if}~\left|\frac{b}{a}\right|\ge b_0,
\end{cases}
$$
where $b_0$ is given in Theorem \ref{thm:main}.
\end{thm}

\def\cprime{$'$} \def\cprime{$'$} \def\cprime{$'$}
\providecommand{\bysame}{\leavevmode\hbox to3em{\hrulefill}\thinspace}
\providecommand{\MR}{\relax\ifhmode\unskip\space\fi MR }
\providecommand{\MRhref}[2]{%
  \href{http://www.ams.org/mathscinet-getitem?mr=#1}{#2}
}
\providecommand{\href}[2]{#2}

\end{document}